\documentclass{IEEEtran}
\usepackage{cite}
\usepackage{amsmath,amssymb,amsfonts}
\usepackage{algorithmic}
\usepackage{graphicx}
\usepackage{subfigure}
\usepackage{textcomp}
\usepackage{xcolor}
\usepackage{siunitx}
\usepackage{todonotes}
\usepackage{comment}
\usepackage{mathtools}
\usepackage{multirow}
\usepackage{multicol}
\usepackage{diagbox}
\usepackage{booktabs}

    \def\BibTeX{{\rm B\kern-.05em{\sc i\kern-.025em b}\kern-.08em
    T\kern-.1667em\lower.7ex\hbox{E}\kern-.125emX}}
    
\begin{document}
\title{On the Computation of Square Roots and Inverse Square Roots of Gram Matrices for Surface Integral Equations in Electromagnetics}
\author{Rui Chen, \IEEEmembership{Member, IEEE}, Adrien Merlini, \IEEEmembership{Senior Member, IEEE}, and Francesco P. Andriulli, \IEEEmembership{Fellow, IEEE}
\thanks{This work has received support by the European Union – Next Generation EU within the PNRR project ``Multiscale modeling and Engineering Applications'' of the Italian National Center for HPC, Big Data and Quantum Computing (Spoke 6) – PNRR M4C2, Investimento 1.4 - Avviso n. 3138 del 16/12/2021 - CN00000013 National Centre for HPC, Big Data and Quantum Computing (HPC) - CUP E13C22000990001, from the European Research Council (ERC) under the European Union’s Horizon 2020 research and innovation programme (grant agreement No. 724846, project 321), from the National Natural Science Foundation of China (NSFC) under Grant 62201264 and Grant 62331016, from Foundation of Anhui Key Laboratory of Industrial Energy-Saving and Safety under Grant KFKT202306, and from the ANR Labex CominLabs under the project ``CYCLE''.  \emph{(Corresponding author: Francesco P. Andriulli.)}}
\thanks{Rui Chen was with the Department of Electronics and Telecommunications, Politecnico di Torino, 10129 Turin, Italy. He is now with the School of Microelectronics (School of Integrated Circuits), Nanjing University of Science and Technology, Nanjing 210094, China and also with North Automatic Control Technology Institute, Taiyuan 030006, China (e-mail: rui.chen@njust.edu.cn). }
\thanks{Adrien Merlini is with 
the Microwave Department, IMT Atlantique, 29238 Brest, France (e-mail: adrien.merlini@imt-atlantique.fr).}
\thanks{Francesco P. Andriulli is with the Department of Electronics and Telecommunications, Politecnico di Torino, 10129 Turin, Italy (e-mail: francesco.andriulli@polito.it). }
}

\maketitle

\begin{abstract}
Surface integral equations (SIEs)-based boundary element methods are widely used for analyzing electromagnetic scattering scenarii. However, after discretization of SIEs, the spectrum and eigenvectors of the boundary element matrices are not usually representative of the spectrum and eigenfunctions of the underlying surface integral operators, which can be problematic for methods that rely heavily on spectral properties. To address this issue,
we delineate some efficient algorithms that allow for the computation of matrix square roots and inverse square roots of the Gram matrices corresponding to the discretization scheme, which can be used for revealing the spectrum of standard electromagnetic integral operators. The algorithms, which are based on properly chosen expansions of the square root and inverse square root functions, are quite effective when applied to
several of the most relevant Gram matrices used
for boundary element discretizations in
electromagnetics.
Tables containing different sets of expansion coefficients are
provided along with  comparative numerical experiments that evidence
advantages and disadvantages of the different approaches. In addition, to
demonstrate the spectrum-revealing properties of the proposed techniques, they are applied to the discretization
of the problem of scattering by a sphere for which the analytic spectrum is known.
\end{abstract}

\begin{IEEEkeywords}
Boundary element method, fast algorithm, Gram matrix, surface integral equation.
\end{IEEEkeywords}

\section{Introduction}
\label{sec:introduction}

In the last few decades, the boundary element method (BEM) has been widely used for analyzing electromagnetic scattering from perfect electrically conducting (PEC) scatterers \cite{Peterson, Gibson, Jin}. For such scatterers, BEM is typically applied to surface integral equations (SIEs), e.g., electric field integral equation (EFIE), magnetic field integral equation, and their linear combinations. In contrast to differential equation-based methods, BEM does not require the use of artificial absorbing boundary conditions for truncating the infinite domain and only requires the discretization of the scatterer's surface with boundary elements, thus resulting in a comparatively lower number of unknowns \cite{Jin}. 

To numerically solve SIEs, surface integral operators are first discretized by expanding the current on the scatterer's surface with basis functions (the boundary elements), then testing functions are selected to test the resulting equation resulting in a fully discretized matrix equation system. Finally,  direct solvers or iterative solvers coupled with fast matrix-vector product algorithms  are used to solve these systems for the unknown current expansion coefficients \cite{Jin}.

For the above-described discretization of SIEs, however, the spectrum and eigenvectors of the boundary element matrices are not usually representative of the spectrum and eigenfunctions of the underlying surface integral operators. This is in part due to the fact that, in general, boundary elements do not form an orthonormal basis of the spaces they span and the operatorial eigenvalues could only be obtained by solving a generalized eigenvalue problem \cite{boffiCompatible}. While this is not necessarily problematic for standard solution strategies, some numerical methods rely on delicate spectral manipulations that necessitate correspondence between the spectra of the operator to be discretized and the corresponding boundary element matrix \cite{Merlini}. 
To address this issue, one could use orthonormal basis functions, but they are often impractical.
Other approaches are instead based on a \textit{a posteriori} orthonormalization of the bases used in the discretization of the integral operators using square roots and inverse square roots of Gram matrices (see \cite{Merlini} and references therein). Several algorithms have been presented to compute these matrices (e.g. \cite{DenmanMatrix,HighamNewton}) and their theory is well-established. Nonetheless, the application of these schemes to the variety of Gram matrices appearing in the integral equations relevant scenarios is far from trivial as a specific set of approaches among many must be selected and the associated parameters need to be optimized to ensure a computational complexity compatible with traditional fast solution and iterative eigenproblem methods \cite{Stewart}.  

In this work, approaches for obtaining square roots and inverse square roots of several of the most relevant Gram matrices are delineated and tailored to all cases of interest. Moreover, the thorough performance analysis is presented together with a comprehensive set of implementation oriented strategies.
To numerically compute the square roots and inverse square roots of the Gram matrices, suitably chosen expansions of the square root and inverse square root scalar functions are used in combination with the theory of matrix functions \cite{Higham}. Several algorithmic approaches for the evaluation of the square roots and inverse square roots of the Gram matrices are considered and tables containing expansion coefficients are provided. Numerical experiments illustrate the merits of the different approaches. Finally, to demonstrate the applicability of the  techniques,  numerical examples showing how the spherical harmonics spectrum of the EFIE operator can be recovered from a standard EFIE system matrix are presented. While most of the discussion focuses on the EFIE and related discretization schemes, the proposed methods are applicable to most of the standard electromagnetic and acoustic integral equations.

\section{Background}
\label{sec:formulation}

Let $S$ denote the surface of a PEC scatterer residing in an unbounded homogeneous background medium. The permittivity and permeability of the background medium are $\varepsilon$ and $\mu$, respectively. The PEC scatterer is excited by an electromagnetic wave with the incident electric field $\mathbf{E}^{\mathrm{inc}}(\mathbf{r})$. Upon excitation, the surface current $\mathbf{J}(\mathbf{r})$ is induced on $S$ that generates the scattered electric field $\mathbf{E}^{\mathrm{sca}}(\mathbf{r})$ in the background medium; $\mathbf{E}^{\mathrm{sca}}(\mathbf{r})$ can be expressed as a function of $\mathbf{J}(\mathbf{r})$ as 
\begin{align}
\nonumber\mathbf{E}^{\mathrm{sca}}(\mathbf{r})=-\mathcal{L}[\mathbf{J}](\mathbf{r})=&-j \omega \mu \int_{S} g (\mathbf{r},\mathbf{r}') \mathbf{J}(\mathbf{r}') ds' \\
&- \frac{j}{\omega \varepsilon} \nabla \int_{S} g (\mathbf{r},\mathbf{r}') \nabla_{\mathrm{s}} ' \cdot \mathbf{J}(\mathbf{r}') ds'.
\label{eq1}
\end{align}
Here, $g(\mathbf{r},\mathbf{r}')=e^{-jkR}/(4\pi R)$ is the Green function, $k=\omega \sqrt{\mu \varepsilon}$ is the wavenumber, and $R=|\mathbf{r}-\mathbf{r}'|$ represents the distance between the field and source points. Note that, the time-harmonic factor $e^{j \omega t}$ is assumed and suppressed throughout this work.

The total electric field $\mathbf{E}^{\mathrm{tot}}(\mathbf{r})=\mathbf{E}^{\mathrm{inc}}(\mathbf{r})+\mathbf{E}^{\mathrm{sca}}(\mathbf{r})$ satisfies the boundary condition $\hat{\mathbf{n}}(\mathbf{r}) \times \mathbf{E}^{\mathrm{tot}}(\mathbf{r}) =0 $ at $\mathbf{r}\in S$, where $\hat{\mathbf{n}}(\mathbf{r})$ is the outward pointing unit normal vector. Inserting \eqref{eq1}  into the above boundary condition yields EFIE as
\begin{align}
\mathcal{T}[\mathbf{J}](\mathbf{r})=-\hat{\mathbf{n}}(\mathbf{r}) \times \mathcal{L}[\mathbf{J}](\mathbf{r})=-\hat{\mathbf{n}}(\mathbf{r}) \times \mathbf{E}^{\mathrm{inc}}(\mathbf{r}).\label{eq2}
\end{align}
To numerically solve \eqref{eq2} for $\mathbf{J}(\mathbf{r})$, $S$ is discretized into a mesh using triangular patches. Then, $\mathbf{J}(\mathbf{r})$ is expanded using Rao-Wilton-Glisson (RWG) basis functions \cite{Rao} as
\begin{align}
\mathbf{J}(\mathbf{r})=\sum_{n=1}^{N_{\mathrm{E}}} \{\mathbf{j}\}_n \mathbf{f}_n(\mathbf{r})\label{eq3}.
\end{align}
Here, $N_{\mathrm{E}}$ denotes the number of edges of the mesh,  $\mathbf{j}$ with dimension $N_{\mathrm{E}} \times 1$ is the vector containing the unknown current expansion coefficients to be solved for, and the $n^{\mathrm{th}}$ RWG basis function $\mathbf{f}_n(\mathbf{r})$ is defined as 
\begin{align}
\mathbf{f}_n(\mathbf{r})=
    \begin{dcases}
      \frac{\mathbf{r}-\mathbf{r}_n^{+}}{2A_n^{+}}, & \mathrm{if} \ \mathbf{r} \in t_n^{+} \\
      \frac{\mathbf{r}_n^{-}-\mathbf{r}}{2A_n^{-}}, & \mathrm{if} \ \mathbf{r} \in t_n^{-}
    \end{dcases},\label{eq4}
\end{align}
where $t_n^{+}$ and $t_n^{-}$ are two triangular patches sharing the $n^{\mathrm{th}}$ edge, $\mathbf{r}_n^{+}$ and $\mathbf{r}_n^{-}$ are the vertices of $t_n^{+}$ and $t_n^{-}$ that do not belong to the $n^{\mathrm{th}}$ edge, and $A_n^{+}$ and $A_n^{-}$ are the areas of $t_n^{+}$ and $t_n^{-}$, respectively.

Next, the rotated RWG functions $\hat{\mathbf{n}}(\mathbf{r}) \times \mathbf{f}_m(\mathbf{r})$, $m=1,2,...,N_{\mathrm{E}}$, are used to test \eqref{eq2}, resulting in the discretized EFIE matrix equation system
\begin{align}
\mathbf{T} \mathbf{j}&=\mathbf{v}.\label{eq5}
\end{align}
Here, the entries of the boundary element matrix $\mathbf{T}$ with dimension $N_{\mathrm{E}} \times N_{\mathrm{E}}$ and the right-hand side vector $\mathbf{v}$ with dimension $N_{\mathrm{E}} \times 1$ are
\begin{align}
\{\mathbf{T}\}_{mn}&=\langle \hat{\mathbf{n}}(\mathbf{r}) \times \mathbf{f}_m(\mathbf{r}), \mathcal{T}[\mathbf{f}_n](\mathbf{r})\rangle\label{eq6}\\
\{\mathbf{v}\}_{m}&=- \langle \hat{\mathbf{n}}(\mathbf{r}) \times \mathbf{f}_m(\mathbf{r}), \hat{\mathbf{n}}(\mathbf{r}) \times \mathbf{E}^{\mathrm{inc}}(\mathbf{r})\rangle\label{eq7},
\end{align}
where $\langle \mathbf{a}(\mathbf{r}), \mathbf{b}(\mathbf{r}) \rangle=\int_{S_{\mathrm{a}}} \mathbf{a}(\mathbf{r}) \cdot \mathbf{b}(\mathbf{r}) ds$ for vector functions, $\langle a(\mathbf{r}), b(\mathbf{r}) \rangle=\int_{S_{\mathrm{a}}} a(\mathbf{r})  b(\mathbf{r}) ds$ for scalar functions and $S_{\mathrm{a}}$ is the support of $\mathbf{a}(\mathbf{r})$ or $a(\mathbf{r})$. 

For certain applications, such as preconditioning \cite{andriulliMultiply}, the electric field integral operator (EFIO) $\mathcal{T}$ must also be discretized using Buffa-Christiansen (BC) basis functions $\mathbf{g}_n(\mathbf{r})$, $n=1,2,...,N_{\mathrm{E}}$ \cite{Buffa}. The BC basis functions are  dual functions with respect to the RWG basis functions and are defined on  barycentrically refined triangular patches. The explicit definition  of  the  BC basis function $\mathbf{g}_n(\mathbf{r})$ is omitted here and can be found in \cite{Buffa}.

We denote the Gram matrices associated with the RWG and BC functions as $\mathbf{G}_{\mathbf{f} , \mathbf{f}}$  and $\mathbf{G}_{\mathbf{g} , \mathbf{g}}$. Their entries are $\{\mathbf{G}_{\mathbf{f} , \mathbf{f}}\}_{mn}=\langle \mathbf{f}_m(\mathbf{r}), \mathbf{f}_n(\mathbf{r})\rangle$ and $\{\mathbf{G}_{\mathbf{g} , \mathbf{g}}\}_{mn}=\langle \mathbf{g}_m(\mathbf{r}), \mathbf{g}_n(\mathbf{r})\rangle$, respectively.

Besides the above mentioned vector basis functions, i.e., the RWG and BC basis functions, used for the discretization of EFIE, scalar basis functions can also be used for the discretization of other SIEs, e.g., potential integral equations  in electromagnetics \cite{Chew, Vico, Rui1} and standard integral equations in acoustics \cite{Rui2, Rui3}. 

In this work, two commonly used scalar basis functions, i.e., pyramid basis functions $\lambda_n(\mathbf{r})$, $n=1,2,...,N_{\mathrm{V}}$ and the dual functions of the pyramid basis functions $\tilde{\lambda}_n(\mathbf{r})$, $n=1,2,...,N_{\mathrm{P}}$, are considered, where  $N_{\mathrm{V}}$ and $N_{\mathrm{P}}$ denote the number of vertices and triangular patches of the mesh, respectively. The dual pyramid basis functions $\tilde{\lambda}_n(\mathbf{r})$ are defined on the barycentrically refined triangular patches. The $n^{\mathrm{th}}$ pyramid basis function $\lambda_n(\mathbf{r})$ is defined to be equal to $1$ at the $n^{\mathrm{th}}$ vertex, $0$ at the other vertices, and linear on the triangular patches sharing the $n^{\mathrm{th}}$ vertex \cite{Chen}. The explicit definition
of the dual pyramid basis function $\tilde{\lambda}_n(\mathbf{r})$ is omitted here that can be found in \cite{Buffa}.
Similarly, the Gram matrices associated with the  pyramid and dual pyramid functions are $\mathbf{G}_{\lambda , \lambda}$ with dimension $N_{\mathrm{V}} \times N_{\mathrm{V}}$ and $\mathbf{G}_{\tilde{\lambda},\tilde{\lambda}}$ with dimension $N_{\mathrm{P}} \times N_{\mathrm{P}}$, whose entries are $\{\mathbf{G}_{\lambda , \lambda}\}_{mn}=\langle \lambda_m(\mathbf{r}), \lambda_n(\mathbf{r})\rangle$ and $\{\mathbf{G}_{\tilde{\lambda} , \tilde{\lambda}}\}_{mn}=\langle \tilde{\lambda}_m(\mathbf{r}), \tilde{\lambda}_n(\mathbf{r})\rangle$, respectively.

The Gram matrices introduced above ($\mathbf{G}_{\mathbf{f} , \mathbf{f}}$, $\mathbf{G}_{\mathbf{g} , \mathbf{g}}$, $\mathbf{G}_{\lambda , \lambda}$, and $\mathbf{G}_{\tilde{\lambda},\tilde{\lambda}}$) are symmetric positive definite (SPD) and can be used for the normalizations of the EFIE boundary element matrix $\mathbf{T}$, the loop-to-RWG transformation matrix $\mathbf{\Lambda}$, and the star-to-RWG transformation matrix $\mathbf{\Sigma}$, where the explicit definitions of $\mathbf{\Lambda}$ and $\mathbf{\Sigma}$ can be found in \cite{SimonInsights}. As mentioned in \cite{Merlini}, using the normalized matrices $\tilde{\mathbf{T}}$, $\tilde{\mathbf{\Lambda}}$, and $\tilde{\mathbf{\Sigma}}$ instead of their standard counterparts $\mathbf{T}$, $\mathbf{\Lambda}$, and $\mathbf{\Sigma}$ allows for more efficient spectral manipulations for problems involving non-uniformly discretized geometries. The normalized matrices are obtained by multiplying the standard matrices with the inverse square roots of the corresponding Gram matrices, e.g., $\tilde{\mathbf{T}} = \sqrt{\mathbf{G}_{\mathbf{f} , \mathbf{f}}}^{-1} \mathbf{T} \sqrt{\mathbf{G}_{\mathbf{f} , \mathbf{f}}}^{-1}$. In Section~\ref{sec:algorithms}, the strategies of the computation of these  square roots and inverse square roots of the Gram matrices are presented.

\section{Computation of the  Square Roots and Inverse Square Roots of the Gram Matrices}
\label{sec:algorithms}

In the following, we will assume to deal with the SPD matrices. In this case, to compute the square roots and inverse square roots of the Gram matrices, a simple-to-define but expensive-to-compute approach is to leverage the singular value decomposition (SVD). The SVD of a real SPD matrix $\mathbf{G}$ of dimension $N \times N$ is defined as 
\begin{align}
\mathbf{G}=\mathbf{U}\mathbf{S}\mathbf{U}^{T},\label{eq8}
\end{align}
where $\mathbf{S}$ is a diagonal matrix of which the elements are the $N$ singular values of $\mathbf{G}$, $\mathbf{U}$ is a unitary matrix, and $\mathbf{G}= \{\mathbf{G}_{\mathbf{f} , \mathbf{f}}, \mathbf{G}_{\mathbf{g} , \mathbf{g}}, \mathbf{G}_{\lambda , \lambda}, \mathbf{G}_{\tilde{\lambda},\tilde{\lambda}} \}$.
Given that the matrices under consideration are also invertible, their square root and inverse square root can be obtained as
\begin{align}
\sqrt{\mathbf{G}}&=\mathbf{U}\bar{\mathbf{S}}\mathbf{U}^{T}\label{eq9}\\
\sqrt{\mathbf{G}}^{-1}&=\mathbf{U}\bar{\mathbf{S}}^{-1}\mathbf{U}^{T},\label{eq10}
\end{align}
where $\bar{\mathbf{S}}$ is a diagonal matrix of which the diagonal entries are the square roots of that of $\mathbf{S}$. However, because of the computational cost of the SVD algorithm, this method quickly becomes unpractical for large problems. Note that, since the matrices under consideration are SPD, their respective square roots are unique \cite{JohnsonUniqueness}.

In the following, the  square roots and inverse square roots of the Gram matrices are expressed as
\begin{align}
\sqrt{\mathbf{G}}&=\sqrt{\|\mathbf{G} \|_2}\sqrt{\mathbf{G}/\|\mathbf{G} \|_2}\,,\label{eq11}\\
\sqrt{\mathbf{G}}^{-1}&=\sqrt{\|\mathbf{G} \|_2}^{-1}\sqrt{\mathbf{G}/{\|\mathbf{G} \|_2}}^{-1}\label{eq12},
\end{align}
where $\| \mathbf{G} \|_2$ denotes the $L_2$-norm of $\mathbf{G}$, which can be computed using power methods. Then, because the Gram matrices under consideration are sparse, $\sqrt{\mathbf{G}/{\|\mathbf{G} \|_2}}$ and $\sqrt{\mathbf{G}/{\|\mathbf{G} \|_2}}^{-1}$ are numerically approximated using polynomial expansions based on the theory of matrix functions \cite{Higham}. Finally, $\sqrt{\mathbf{G}}$ and $\sqrt{\mathbf{G}}^{-1}$ are computed by multiplying  $\sqrt{\|\mathbf{G} \|_2}$ and $\sqrt{\|\mathbf{G} \|_2}^{-1}$ with the numerical approximations of $\sqrt{\mathbf{G}/{\|\mathbf{G} \|_2}}$ and $\sqrt{\mathbf{G}/{\|\mathbf{G} \|_2}}^{-1}$, respectively.

To further condense the  treatments, we denote a matrix function by a single symbol $f$ that  corresponds to the matrix function either $\mathbf{X} \mapsto \sqrt{\mathbf{X}}$ or $\mathbf{X} \mapsto \sqrt{\mathbf{X}}^{-1}$ defined on the space of the SPD matrices of dimension $N \times N$. The numerical expression of $f(\mathbf{X})$ can be obtained by generalizing the numerical (polynomial or rational) approximation of its scalar function counterpart either $x \mapsto \sqrt{x}$ or $x \mapsto \sqrt{x}^{-1}$, $x\in\mathbb{R}^+$, respectively, to the matrix argument \cite{Higham}. In the next Section~\ref{subsec:TSE}, \ref{subsec:CPE}, and \ref{subsec:PAE}, three algorithmic approaches based on the polynomial approximations are presented to achieve this purpose.

\subsection{Taylor Series Expansion}
\label{subsec:TSE}

As the first algorithmic approach, the Taylor series expansion (TSE) \cite{Taylor} is used for the polynomial approximations of $x \mapsto \sqrt{x}$ and $x \mapsto \sqrt{x}^{-1}$. Using TSE, $x \mapsto \sqrt{x}$ and $x \mapsto \sqrt{x}^{-1}$ can be numerically expanded around $x=1$ as 
\begin{align}
x \mapsto
    \begin{dcases}
      \sqrt{x}\approx\sum_{n=0}^{N_{\mathrm{T}}} c^{\mathrm{TSE}}_n (x-1)^n \\
      \sqrt{x}^{-1}\approx\sum_{n=0}^{N_{\mathrm{T}}} \tilde{c}^{\mathrm{TSE}}_n (x-1)^n
    \end{dcases}.\label{eq13}
\end{align}
Here, $N_{\mathrm{T}}$ is the TSE expansion order, and $c^{\mathrm{TSE}}_n$ and $\tilde{c}^{\mathrm{TSE}}_n$ are the $n^{\mathrm{th}}$ expansion coefficients of TSE for the approximations of $\sqrt{x}$ and $\sqrt{x}^{-1}$, respectively, which can be computed as
\begin{equation}
c^{\mathrm{TSE}}_n=\frac{\sqrt{x}^{(n)}|_{x=1}}{n!}\quad \text{and} \quad \tilde{c}^{\mathrm{TSE}}_n =\frac{(1/\sqrt{x})^{(n)}|_{x=1}}{n!}.\label{eq14}
\end{equation}
The approximations of their matrix function counterparts $\mathbf{X} \mapsto \sqrt{\mathbf{X}}$ and $\mathbf{X} \mapsto \sqrt{\mathbf{X}}^{-1}$ are obtained by generalizing \eqref{eq13} to the matrix arguments as
\begin{equation}
\mathbf{X} \mapsto
    \begin{dcases}
      \sqrt{\mathbf{X}}\approx\sum_{n=0}^{N_{\mathrm{T}}} c^{\mathrm{TSE}}_n (\mathbf{X}-\mathbf{I})^n \\
      \sqrt{\mathbf{X}}^{-1}\approx\sum_{n=0}^{N_{\mathrm{T}}} \tilde{c}^{\mathrm{TSE}}_n (\mathbf{X}-\mathbf{I})^n
    \end{dcases}.\label{eq15}
\end{equation}
Here, $\mathbf{I}$ is the identity matrix with the same dimension as $\mathbf{X}$. The values of the first ten expansion coefficients of TSE $\{c^{\mathrm{TSE}}_0,c^{\mathrm{TSE}}_1,c^{\mathrm{TSE}}_2,c^{\mathrm{TSE}}_3,c^{\mathrm{TSE}}_4,c^{\mathrm{TSE}}_5,c^{\mathrm{TSE}}_6,c^{\mathrm{TSE}}_7,c^{\mathrm{TSE}}_8,c^{\mathrm{TSE}}_9\}$ and $\{\tilde{c}^{\mathrm{TSE}}_0,\tilde{c}^{\mathrm{TSE}}_1,\tilde{c}^{\mathrm{TSE}}_2,\tilde{c}^{\mathrm{TSE}}_3,\tilde{c}^{\mathrm{TSE}}_4,\tilde{c}^{\mathrm{TSE}}_5,\tilde{c}^{\mathrm{TSE}}_6,\tilde{c}^{\mathrm{TSE}}_7,\tilde{c}^{\mathrm{TSE}}_8,\tilde{c}^{\mathrm{TSE}}_9\}$ for the approximations of $\sqrt{\mathbf{X}}$ and $\sqrt{\mathbf{X}}^{-1}$ are $\{1, \frac{1}{2}, -\frac{1}{8}, \frac{1}{16}, -\frac{5}{128}, \frac{7}{256}, -\frac{21}{1024}, \frac{33}{2048}, -\frac{429}{32768}, \frac{715}{65536}\}$ and $\{1, -\frac{1}{2}, \frac{3}{8}, -\frac{5}{16}, \frac{35}{128}, -\frac{63}{256}, \frac{231}{1024}, -\frac{429}{2048}, \frac{6435}{32768}, -\frac{12155}{65536}\}$, respectively.

\subsection{Chebyshev Polynomial Expansion}
\label{subsec:CPE}

As the second algorithmic approach, the truncated Chebyshev polynomial series expansion (CPE) \cite{Mason} is used for the polynomial approximations of the square root function and its reciprocal on the interval $x \in [n_0,1]$ as
\begin{align}
x \mapsto
    \begin{dcases}
      \sqrt{x}\approx\sum_{n=0}^{N_{\mathrm{C}}}{}^\prime c^{\mathrm{CPE}}_n T_n(x) \\
      \sqrt{x}^{-1}\approx\sum_{n=0}^{N_{\mathrm{C}}}{}^\prime \tilde{c}^{\mathrm{CPE}}_n T_n(x)
    \end{dcases}.\label{eq16}
\end{align}
Here, the primed summation indicates that the first term is halved, $n_0$ is a positive real number, $N_{\mathrm{C}}$ is the expansion order of CPE, 
and $T_n(x)$  is the shifted Chebyshev polynomial of the first kind of order $n$ defined on the interval $x \in [n_0,1]$  as
\begin{align}
T_n(x)&=
    \begin{dcases}
      1 & \mathrm{if} \ n=0 \\
      \frac{2x-(n_0+1)}{1-n_0} & \mathrm{if} \ n=1 \\
      2\left(\frac{2x-(n_0+1)}{1-n_0}\right)T_{n-1}(x)-T_{n-2}(x) & \mathrm{otherwise}
    \end{dcases}.\label{eq17}
\end{align}
The Chebyshev series coefficients $c^{\mathrm{CPE}}_n$ and $\tilde{c}^{\mathrm{CPE}}_n$ for the approximations of $\sqrt{x}$ and $\sqrt{x}^{-1}$ can be computed as \cite{Mason}
\begin{align}
c^{\mathrm{CPE}}_n&=\frac{2}{\pi} \int_{n_0}^1 \frac{\sqrt{x} T_n(x)}{\sqrt{-x^2+(n_0+1)x-n_0}}\mathrm{d}x\label{eq18}\,\\
\tilde{c}^{\mathrm{CPE}}_n&=\frac{2}{\pi} \int_{n_0}^1 \frac{\sqrt{x}^{-1} T_n(x)}{\sqrt{-x^2+(n_0+1)x-n_0}}\mathrm{d}x\,,\label{eq19}
\end{align}
which can be rewritten with special functions as
\begin{align}
c^{\mathrm{CPE}}_n &= 2 \sqrt{n_0}  \frac{\,{}_3F_2\left(\begin{matrix}\begin{matrix}-1/2 & 1/2 & 1\end{matrix}\\\begin{matrix}1 - n & 1 + n\end{matrix}\end{matrix}; 1 - n_0^{-1}\right)}{\Gamma(1 - n)\Gamma(1 + n)}\label{eq:coefs_cheby_sqrt}\,\\
\tilde{c}^{\mathrm{CPE}}_n &= \frac{2}{\sqrt{n_0}}  \frac{\,{}_3F_2\left(\begin{matrix}\begin{matrix}1/2 & 1/2 & 1\end{matrix}\\\begin{matrix}1 - n & 1 + n\end{matrix}\end{matrix}; 1 - n_0^{-1}\right)}{\Gamma(1 - n)\Gamma(1 + n)} \,,\label{eq:coefs_cheby_isqrt}
\end{align}
where ${}_pF_q$ is the generalized hypergeometric function \cite[eq. (16.2.1)]{NIST} and $\Gamma$ is the Euler gamma function.
 
Similarly, after the series expansions of $x \mapsto \sqrt{x}$ and $x \mapsto \sqrt{x}^{-1}$ in \eqref{eq16} are derived, their matrix function counterparts can be obtained by generalizing \eqref{eq16} to the matrix arguments as
\begin{align}
\mathbf{X} \mapsto
    \begin{dcases}
      \sqrt{\mathbf{X}}\approx\sum_{n=0}^{N_{\mathrm{C}}}{}' c^{\mathrm{CPE}}_n T_n(\mathbf{X}) \\
      \sqrt{\mathbf{X}}^{-1}\approx\sum_{n=0}^{N_{\mathrm{C}}}{}' \tilde{c}^{\mathrm{CPE}}_n T_n(\mathbf{X})
    \end{dcases},\label{eq22}
\end{align}
where $T_n(\mathbf{X})$ is defined as
\begin{align}
T_n(\mathbf{X})&=
    \begin{dcases}
      \mathbf{I} & \mathrm{if} \ n=0 \\
      \frac{2\mathbf{X}-(n_0+1)\mathbf{I}}{1-n_0} & \mathrm{if} \ n=1 \\
      \frac{2[2\mathbf{X}-(n_0+1)\mathbf{I}]}{1-n_0}T_{n-1}(\mathbf{X})-T_{n-2}(\mathbf{X}) & \mathrm{otherwise}
    \end{dcases}.\label{eq23}
\end{align}
Given \eqref{eq11} and \eqref{eq12}, the spectrum of the SPD matrix to which CPE is applied lives in the interval $[n_0, 1]$, where $n_0$ is the inverse of the condition number of the matrix.

When using \eqref{eq22}, the coefficients $c^{\mathrm{CPE}}_n$ and $\tilde{c}^{\mathrm{CPE}}_n$ must be computed for each matrix, since they are dependent on the condition number $n_0^{-1}$ of the matrix argument. If this turns out to be inconvenient, tabulated values of coefficients can be used for the matrices that have a value of $n_0$ living in known ranges. Usage of such tabulated values will typically yield higher errors than the computation with the matrix-dependent coefficients and in Section \ref{sec:results} we will provide numerical results for these two approaches designated respectively as ``CPE-1'' for the matrix specific approach and ``CPE-2'' for the tabulated approach. For convenience, the values of $c^{\mathrm{CPE}}_n$ and $\tilde{c}^{\mathrm{CPE}}_n$ of CPE-2 are provided in Table~\ref{tab:coefs_chebyshev_sqrt} and Table~\ref{tab:coefs_chebyshev_isqrt} for given bounds on the Gram matrices conditioning (different $n_0$).
The order $N_{\mathrm{C}}$ at which the Chebyshev series should be truncated to reach a target relative error $\delta$ on the matrix (defined in Section \ref{sec:results}) has been obtained through numerical analyses and is reported in Table~\ref{tab:order_chebyshev_sqrt}. Note that, the direct evaluation of the sum in \eqref{eq22} is not necessarily optimal and other algorithms can be employed \cite{clenshaw1955note}.

\subsection{Pad\'e Approximant Expansion}
\label{subsec:PAE}

The last algorithmic approach we explore 
is the Pad\'e approximant expansion (PAE) \cite{Baker} for $x \mapsto \sqrt{x}$ and $x \mapsto \sqrt{x}^{-1}$ at $x=1$ as
\begin{align}
x \mapsto
    \begin{dcases}
      \sqrt{x}\approx\dfrac{\sum_{n=0}^{N_{\mathrm{A}}} c^{\mathrm{PAE}}_n x^n}{\sum_{n=0}^{N_{\mathrm{A}}} c^{\mathrm{PAE}}_n x^{N_{\mathrm{A}}-n}} \\
      \sqrt{x}^{-1}\approx\dfrac{\sum_{n=0}^{N_{\mathrm{A}}} c^{\mathrm{PAE}}_n x^{N_{\mathrm{A}}-n}}{\sum_{n=0}^{N_{\mathrm{A}}} c^{\mathrm{PAE}}_n x^n}
    \end{dcases}.\label{eq24}
\end{align}
Here, $N_{\mathrm{A}}$ is the expansion order of PAE, and $c^{\mathrm{PAE}}_n$ is the $n^{\mathrm{th}}$ expansion coefficient of PAE for the approximations of $\sqrt{x}$ and $\sqrt{x}^{-1}$, which depends on $N_{\mathrm{A}}$ and can be obtained after solving the linear equation systems as presented in \cite{Song}. The numerical expressions of their matrix function counterparts are obtained by generalizing \eqref{eq24} to the matrix arguments as
\begin{align}
\mathbf{X} \mapsto
    \begin{dcases}
      \sqrt{\mathbf{X}}\approx \left(\sum_{n=0}^{N_{\mathrm{A}}} c^{\mathrm{PAE}}_n   \mathbf{X}^{N_{\mathrm{A}}-n}\right)^{-1} \left(\sum_{n=0}^{N_{\mathrm{A}}} c^{\mathrm{PAE}}_n \mathbf{X}^n\right) \\
      \sqrt{\mathbf{X}}^{-1}\approx \left(\sum_{n=0}^{N_{\mathrm{A}}} c^{\mathrm{PAE}}_n \mathbf{X}^n\right)^{-1} \left(\sum_{n=0}^{N_{\mathrm{A}}} c^{\mathrm{PAE}}_n   \mathbf{X}^{N_{\mathrm{A}}-n}\right)
    \end{dcases}.\label{eq25}
\end{align}
The detailed values of the expansion coefficients $c^{\mathrm{PAE}}_n$ of PAE  for $N_\mathrm{A}=0,1,\dots,9$ are provided in Table~\ref{tab1}.

\begin{table*}
\centering
\caption{First Chebyshev series coefficients for the square root function for different $n_0$.}
\label{tab:coefs_chebyshev_sqrt}
\resizebox{\textwidth}{!}
{
\begin{tabular}{|c|c| c |c |c|c| c |c |c|c |c|}
  \toprule
  \multirow{2}{*}{$n_0$} & $c^{\mathrm{CPE}}_0$ & $c^{\mathrm{CPE}}_1$ & $c^{\mathrm{CPE}}_2$ & $c^{\mathrm{CPE}}_3$ & $c^{\mathrm{CPE}}_4$ & $c^{\mathrm{CPE}}_5$ & $c^{\mathrm{CPE}}_6$ & $c^{\mathrm{CPE}}_7$ & $c^{\mathrm{CPE}}_8$ & $c^{\mathrm{CPE}}_9$ \\[1ex]
  & $c^{\mathrm{CPE}}_{10}$ & $c^{\mathrm{CPE}}_{11}$ & $c^{\mathrm{CPE}}_{12}$ & $c^{\mathrm{CPE}}_{13}$ & $c^{\mathrm{CPE}}_{14}$ & $c^{\mathrm{CPE}}_{15}$ & $c^{\mathrm{CPE}}_{16}$ & $c^{\mathrm{CPE}}_{17}$ & $c^{\mathrm{CPE}}_{18}$ & $c^{\mathrm{CPE}}_{19}$\\
 \midrule
 \midrule
  \multirow{2}{*}{$n_0 \geq \num{1e-1}$} &$\frac{50720584}{36057897}$&$\frac{58192354}{176379781}$&$-\frac{4460738}{108096729}$&$\frac{599673}{57061255}$&$-\frac{287819}{85384633}$&$\frac{183411}{150871267}$&$-\frac{58340}{123912171}$&$\frac{51981}{271721903}$&$-\frac{7566}{94035965}$&$\frac{2764}{79582277}$\\[1ex]
  & $-\frac{201}{13137223}$&$\frac{1188}{173409911}$&$-\frac{409}{131560083}$&$\frac{145}{101639721}$&$-\frac{17}{25723455}$&$\frac{55}{178196989}$&$-\frac{15}{103329721}$&$\frac{3}{43668719}$&$-\frac{1}{30591204}$&$\frac{1}{63979016}$\\
  \midrule
  \multirow{2}{*}{$n_0 \geq \num{5e-2}$} & $\frac{141760564}{104989387}$&$\frac{35002745}{95238932}$&$-\frac{5385084}{98048633}$&$\frac{375119}{22232916}$&$-\frac{123075}{18770531}$&$\frac{964144}{335715223}$&$-\frac{99240}{73358167}$&$\frac{61845}{92416364}$&$-\frac{49497}{144355439}$&$\frac{43265}{239826181}$\\[1ex]
  &$-\frac{8116}{83759003}$ & $\frac{4645}{87782356}$&$-\frac{4539}{154958201}$&$\frac{2159}{131648930}$&$-\frac{567}{61162961}$&$\frac{708}{133999337}$&$-\frac{225}{74184637}$&$\frac{187}{106738132}$&$-\frac{91}{89427045}$&$\frac{77}{129639190}$\\
  \midrule
  \multirow{2}{*}{$n_0 \geq \num{1e-2}$} & $\frac{94293998}{72892523}$&$\frac{11654085}{28548136}$&$-\frac{13024699}{174927782}$&$\frac{1616566}{56730601}$&$-\frac{1147126}{82276417}$&$\frac{761833}{98466580}$&$-\frac{308484}{66539819}$&$\frac{263708}{90144229}$&$-\frac{162889}{85022700}$&$\frac{161374}{125071371}$\\[1ex]
  &$-\frac{37949}{42730513}$&$\frac{43604}{70096297}$&$-\frac{29207}{66081562}$&$\frac{51537}{162168061}$&$-\frac{25435}{110193549}$&$\frac{11535}{68213866}$&$-\frac{10239}{82034690}$&$\frac{3604}{38866185}$&$-\frac{6215}{89696092}$&$\frac{2753}{52900755}$\\
  \midrule
  \multirow{2}{*}{$n_0 \geq \num{5e-3}$} & $\frac{111282682}{86634283}$&$\frac{19778972}{47632141}$&$-\frac{10097893}{128419878}$&$\frac{1613988}{51323195}$&$-\frac{1103585}{68372228}$&$\frac{1156783}{122811139}$&$-\frac{1072665}{180416284}$&$\frac{182246}{46057955}$&$-\frac{335564}{122664125}$&$\frac{102602}{52715219}$\\[1ex]
  &$-\frac{86227}{60890409}$&$\frac{77118}{73520351}$&$-\frac{61735}{78300108}$&$\frac{43125}{71884333}$&$-\frac{93187}{202045214}$&$\frac{27405}{76607114}$&$-\frac{26887}{96161515}$&$\frac{1853}{8422610}$&$-\frac{19053}{109416983}$&$\frac{2101}{15164315}$\\
  \midrule
  \multirow{2}{*}{$n_0 \geq \num{1e-3}$} & $\frac{67448647}{52859297}$&$\frac{18465778}{43750239}$&$-\frac{3199979}{38492102}$&$\frac{3836927}{109928911}$&$-\frac{1863218}{98477081}$&$\frac{1438469}{122838561}$&$-\frac{1071366}{136237541}$&$\frac{1846182}{330775465}$&$-\frac{1551939}{376426982}$&$\frac{400555}{127615892}$\\[1ex]
  &$-\frac{362239}{148056512}$&$\frac{272353}{140130127}$&$-\frac{353347}{225334648}$&$\frac{98124}{76560973}$&$-\frac{75049}{70863879}$&$\frac{78769}{89167398}$&$-\frac{41954}{56475965}$&$\frac{258049}{410152123}$&$-\frac{33661}{62777131}$&$\frac{177773}{386862333}$\\
  \bottomrule
\end{tabular}
}
\end{table*}

\begin{table*}
\centering
\caption{First Chebyshev series coefficients for the reciprocal square root function for different $n_0$.}
\label{tab:coefs_chebyshev_isqrt}
\resizebox{\textwidth}{!}
{
\begin{tabular}{|c|c| c |c |c|c| c |c |c|c |c|}
  \toprule
  \multirow{2}{*}{$n_0$} & $\tilde{c}^{\mathrm{CPE}}_0$ & $\tilde{c}^{\mathrm{CPE}}_1$ & $\tilde{c}^{\mathrm{CPE}}_2$ & $\tilde{c}^{\mathrm{CPE}}_3$ & $\tilde{c}^{\mathrm{CPE}}_4$ & $\tilde{c}^{\mathrm{CPE}}_5$ & $\tilde{c}^{\mathrm{CPE}}_6$ & $\tilde{c}^{\mathrm{CPE}}_7$ & $\tilde{c}^{\mathrm{CPE}}_8$ & $\tilde{c}^{\mathrm{CPE}}_9$ \\[1ex]
  & $\tilde{c}^{\mathrm{CPE}}_{10}$ & $\tilde{c}^{\mathrm{CPE}}_{11}$ & $\tilde{c}^{\mathrm{CPE}}_{12}$ & $\tilde{c}^{\mathrm{CPE}}_{13}$ & $\tilde{c}^{\mathrm{CPE}}_{14}$ & $\tilde{c}^{\mathrm{CPE}}_{15}$ & $\tilde{c}^{\mathrm{CPE}}_{16}$ & $\tilde{c}^{\mathrm{CPE}}_{17}$ & $\tilde{c}^{\mathrm{CPE}}_{18}$ & $\tilde{c}^{\mathrm{CPE}}_{19}$\\
 \midrule
 \midrule
 \multirow{2}{*}{$n_0 \geq \num{1e-1}$} &$\frac{374048017}{113951175}$&$-\frac{88499941}{99875048}$&$\frac{15797788}{45155977}$&$-\frac{20535446}{134671187}$&$\frac{3198561}{45955201}$&$-\frac{1821274}{55810317}$&$\frac{1646349}{105728006}$&$-\frac{656367}{87247730}$&$\frac{187127}{51011883}$&$-\frac{257569}{142972536}$\\[1ex]
 &$\frac{337750}{379582213}$&$-\frac{17881}{40498318}$&$\frac{23207}{105516984}$&$-\frac{16946}{154174661}$&$\frac{5735}{104114243}$&$-\frac{1684}{60855691}$&$\frac{1403}{100712358}$&$-\frac{640}{91088439}$&$\frac{138}{38877931}$&$-\frac{197}{109695708}$\\
\midrule
\multirow{2}{*}{$n_0 \geq \num{5e-2}$} &$\frac{249779459}{67453088}$&$-\frac{350792999}{280590374}$&$\frac{132037463}{217146298}$&$-\frac{17309909}{53231083}$&$\frac{24420789}{134319148}$&$-\frac{11726340}{112404257}$&$\frac{9883564}{162320651}$&$-\frac{3377515}{93897642}$&$\frac{1014741}{47325860}$&$-\frac{16780163}{1303748560}$\\[1ex]
&$\frac{510107}{65658935}$&$-\frac{491723}{104378654}$&$\frac{360646}{125772587}$&$-\frac{240771}{137509316}$&$\frac{128666}{120011707}$&$-\frac{79363}{120608489}$&$\frac{30997}{76590848}$&$-\frac{26019}{104339573}$&$\frac{15878}{103168691}$&$-\frac{17123}{180008350}$\\
\midrule
\multirow{2}{*}{$n_0 \geq \num{1e-2}$} &$\frac{237825977}{50542861}$&$-\frac{39264271}{17952244}$&$\frac{195265677}{138817343}$&$-\frac{103416965}{105120298}$&$\frac{12141478}{16961213}$&$-\frac{81979913}{153769588}$&$\frac{100102199}{248249039}$&$-\frac{41303408}{133947555}$&$\frac{47672400}{200509969}$&$-\frac{10291450}{55779169}$\\[1ex]
&$\frac{11305153}{78552187}$&$-\frac{9467215}{83975554}$&$\frac{3834427}{43265789}$&$-\frac{7690733}{110058353}$&$\frac{5948375}{107682632}$&$-\frac{3503490}{80051251}$&$\frac{1048651}{30183225}$&$-\frac{6557825}{237360333}$&$\frac{19174519}{871395101}$&$-\frac{999069}{56928095}$\\
\midrule
\multirow{2}{*}{$n_0 \geq \num{5e-3}$} &$\frac{564146905}{109693414}$&$-\frac{271445188}{103894561}$&$\frac{203729517}{112913626}$&$-\frac{84437949}{62627027}$&$\frac{49111051}{46961990}$&$-\frac{15609456}{18825521}$&$\frac{63222034}{94771407}$&$-\frac{56599573}{104360532}$&$\frac{38412599}{86436742}$&$-\frac{62114496}{169532831}$\\[1ex]
&$\frac{75488883}{248678426}$&$-\frac{2511703}{9946203}$&$\frac{16184137}{76778415}$&$-\frac{21196048}{120119399}$&$\frac{17846601}{120516115}$&$-\frac{7550133}{60623138}$&$\frac{11163511}{106379007}$&$-\frac{5209630}{58817269}$&$\frac{7231330}{96584953}$&$-\frac{11584499}{182801188}$\\
\midrule
\multirow{2}{*}{$n_0 \geq \num{1e-3}$} &$\frac{630048624}{102215551}$&$-\frac{343367843}{94808470}$&$\frac{294957839}{105949003}$&$-\frac{494742647}{216022718}$&$\frac{1224960907}{629659811}$&$-\frac{206423487}{122566312}$&$\frac{180757935}{122418742}$&$-\frac{130317757}{99758659}$&$\frac{119757642}{102912463}$&$-\frac{36868111}{35375072}$\\[1ex]
&$\frac{65195223}{69543947}$&$-\frac{95599410}{112964291}$&$\frac{184341026}{240571663}$&$-\frac{55931467}{80408470}$&$\frac{68082220}{107582399}$&$-\frac{70019489}{121380714}$&$\frac{69973473}{132846113}$&$-\frac{35406207}{73506184}$&$\frac{53675438}{121692277}$&$-\frac{34220975}{84624381}$\\
 \bottomrule
\end{tabular}
}
\end{table*}

\begin{table*}
\centering
\caption{Chebyshev series truncation order for the square root (left) and reciprocal square root (right) functions for different $n_0$ and target matrix relative error $\delta$.}
\label{tab:order_chebyshev_sqrt}
\begin{tabular}{|c|c|c|c|c|c|}
\toprule
\diagbox{$n_0$}{$\delta$} & $10^{-2}$ & $10^{-3}$ & $10^{-4}$ & $10^{-5}$ & $10^{-6}$ \\
\midrule
$n_0 \geq \num{1e-1}$ & 3 & 5  & 8 & 11 & 14 \\
$n_0 \geq \num{5e-2}$ & 4 & 7 & 11 & 15 & 19 \\
$n_0 \geq \num{1e-2}$ & 6 & 13 & 21 & 31 & 40 \\
$n_0 \geq \num{5e-3}$ & 8 & 17 & 29 & 41 & ----- \\
$n_0 \geq \num{1e-3}$ & 12 & 30 & ----- & ----- & ----- \\
\bottomrule
\end{tabular}
~
\begin{tabular}{|c|c|c|c|c|c|}
\toprule
\diagbox{$n_0$}{$\delta$} & $10^{-2}$ & $10^{-3}$ & $10^{-4}$ & $10^{-5}$ & $10^{-6}$ \\
\midrule
$n_0 \geq \num{1e-1}$ & 5 & 9 & 12 & 15 & 19 \\
$n_0 \geq \num{5e-2}$ & 8 & 13 & 17 & 22 & 27 \\
$n_0 \geq \num{1e-2}$ & 18 & 28 & 39 & ----- & ----- \\
$n_0 \geq \num{5e-3}$ & 25 & 40 & ----- & ----- & ----- \\
$n_0 \geq \num{1e-3}$ & ----- & ----- & ----- & ----- & ----- \\
\bottomrule
\end{tabular}
\end{table*}

\begin{table*}
\centering
\caption{Pad\'e approximant expansion coefficients for $N_\mathrm{A}=0,1,\dots,9$.}
\resizebox{0.45\textwidth}{!}{\begin{tabular}{|c| c |c |c|c| c |c |c|c |c|c|} 
 \toprule
  $N_\mathrm{A}$ & $0$ & $1$ & $2$ & $3$ & $4$ & $5$ & $6$ & $7$ & $8$ & $9$ \\  
 \midrule
   $c^{\mathrm{PAE}}_0$ & $1$ & $1$ & $1$ & $1$ & $1$ & $1$ & $1$ & $1$ & $1$ & $1$ \\ 
 $c^{\mathrm{PAE}}_1$ & ----- & $3$ & $10$ & $21$ & $36$ & $55$ & $78$ & $105$ & $136$ & $171$ \\ 
 $c^{\mathrm{PAE}}_2$ & ----- & ----- & $5$ &$35$  & $126$ & $330$ & $715$ & $1365$ & $2380$ & $3876$ \\ 
 $c^{\mathrm{PAE}}_3$ & ----- & ----- & ----- & $7$ &$84$  & $462$ & $1716$ & $5005$ & $12376$ & $27132$ \\ 
 $c^{\mathrm{PAE}}_4$ & ----- & ----- & ----- & ----- & $9$ &$165$  &$1287$  & $6435$ & $24310$ & $75582$ \\ 
 $c^{\mathrm{PAE}}_5$ & ----- & ----- & ----- & ----- & ----- & $11$ & $286$ & $3003$ & $19448$ &  $92378$\\ 
 $c^{\mathrm{PAE}}_6$ & ----- & ----- & ----- & ----- & ----- & ----- & $13$ & $455$ & $6188$ & $50388$ \\ 
 $c^{\mathrm{PAE}}_7$ & ----- & ----- & ----- & ----- & ----- & ----- & ----- & $15$ & $680$ & $11628$ \\ 
 $c^{\mathrm{PAE}}_8$ & ----- &  ----- & ----- & ----- & ----- & ----- & ----- & ----- & $17$ & $969$ \\ 
 $c^{\mathrm{PAE}}_9$ & ----- &-----  & ----- & ----- & ----- & ----- & ----- & ----- & ----- & $19$ \\ 
 \bottomrule
\end{tabular}}
\label{tab1}
\end{table*}

\subsection{General Comments}
\label{subsec:comments}

Before the presentation of the numerical experiments, we would like to attract the reader's attention on a few points.
After the numerical expressions of $\sqrt{\mathbf{X}}$ and $\sqrt{\mathbf{X}}^{-1}$ in \eqref{eq15}, \eqref{eq22}, and \eqref{eq25}  are obtained, the numerical expressions of $\sqrt{\mathbf{G}/{\|\mathbf{G} \|_2}}$ and $\sqrt{\mathbf{G}/{\|\mathbf{G} \|_2}}^{-1}$ are obtained 
 with the substitution $\mathbf{X}=\mathbf{G}/{\|\mathbf{G} \|_2}$   in  \eqref{eq15}, \eqref{eq22}, and \eqref{eq25} for TSE, CPE, and PAE, respectively.
In addition, the normalization of $\mathbf{G}$ (i.e., dividing $\mathbf{G}$ by $\|\mathbf{G} \|_2$) in \eqref{eq11} and \eqref{eq12} is not a necessary step of our proposed algorithms.
However, this normalization is often convenient for implementation and convergence purposes; when the matrices of interest are sparse, SPD, and well-conditioned, the power methods provide a very efficient way for computing $\|\mathbf{G} \|_2$.

\section{Numerical Results}
\label{sec:results}

In this section, several numerical experiments are presented to investigate the advantages and disadvantages of the proposed algorithmic approaches for the computation of $\sqrt{\mathbf{G}}$ and $\sqrt{\mathbf{G}}^{-1}$, $\mathbf{G} \in \{\mathbf{G}_{\mathbf{f} , \mathbf{f}}, \mathbf{G}_{\mathbf{g} , \mathbf{g}}, \mathbf{G}_{\lambda , \lambda}, \mathbf{G}_{\tilde{\lambda},\tilde{\lambda}} \}$. To investigate the accuracy of the numerical approximations of the square roots and inverse square roots of the Gram matrices, the relative error $\delta$ is defined as
\begin{align}
\delta=\frac{\|f^{\mathrm{num}}(\mathbf{G})-f^{\mathrm{ref}}(\mathbf{G}) \|_2}{\|f^{\mathrm{ref}}(\mathbf{G}) \|_2}.\label{eq26}
\end{align}
Here, $f$ denotes $\mathbf{G} \mapsto \sqrt{\mathbf{G}}$ or $\mathbf{G} \mapsto \sqrt{\mathbf{G}}^{-1}$, and the superscripts ``num'' and ``ref'' represent the numerical and reference values of $f(\mathbf{G})$  computed using the proposed   algorithmic approaches in Section~\ref{sec:algorithms} and the computationally intensive method via SVD in \eqref{eq9} and \eqref{eq10}, respectively.

For all the numerical experiments, a PEC sphere of radius \SI{0.5}{\meter}, centered at the origin, and residing in  free space is considered. The surface of the PEC sphere is non-uniformly discretized into three sets of meshes, denoted by ``mesh-1'', ``mesh-2'', and ``mesh-3'' with  $\{ N_{\mathrm{E}}, N_{\mathrm{P}}, N_{\mathrm{V}}\}=\{ 1230, 820, 412  \}$, $\{ N_{\mathrm{E}}, N_{\mathrm{P}}, N_{\mathrm{V}}\}=\{ 2193, 1462, 733  \}$, and $\{ N_{\mathrm{E}}, N_{\mathrm{P}}, N_{\mathrm{V}}\}=\{ 4449, 2966, 1485\}$, respectively. The three sets of meshes have been selected so that the resulting Gram matrices are increasingly ill-conditioned for mesh-1, mesh-2, and mesh-3 with the condition numbers of $\{\mathbf{G}_{\mathbf{f} , \mathbf{f}}, \mathbf{G}_{\mathbf{g} , \mathbf{g}}, \mathbf{G}_{\lambda , \lambda}, \mathbf{G}_{\tilde{\lambda},\tilde{\lambda}}\}$ as $\{5.16,5.74,6.72,9.11 \}$, $\{20.45,14.49,29.12,24.58 \}$, and $\{615.61,344.95,361.91,269.31 \}$, respectively. 

\subsection{Accuracy of the Computation of the  Square Roots and Inverse  Square Roots of the Gram Matrices}

\begin{figure}
\centering
\subfigure[]{\includegraphics[width=0.9\columnwidth]{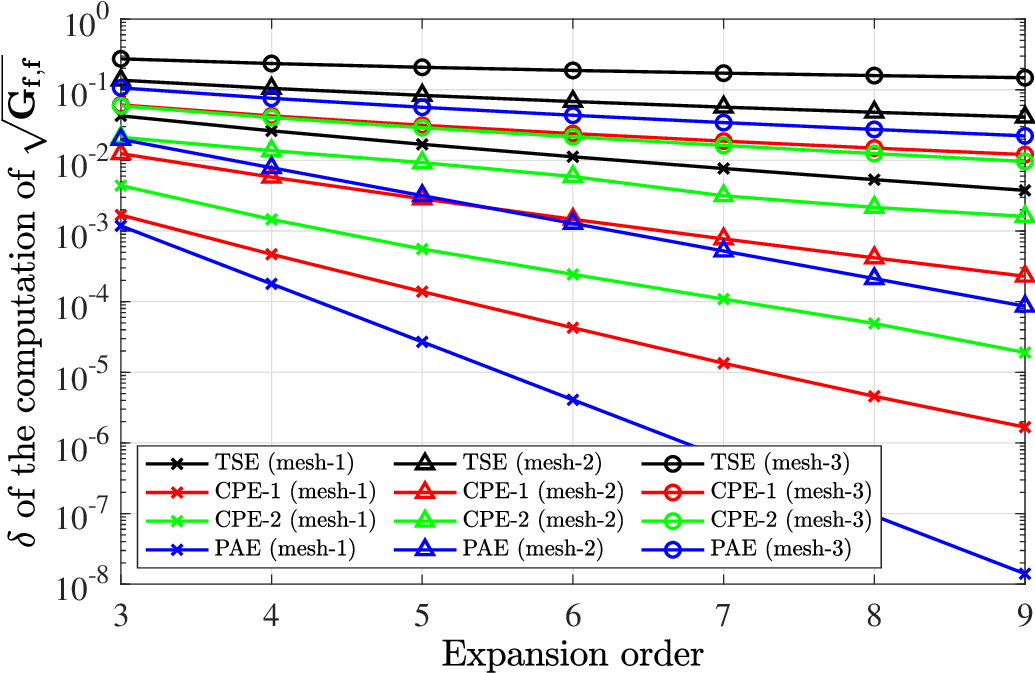}}
\subfigure[]{\includegraphics[width=0.9\columnwidth]{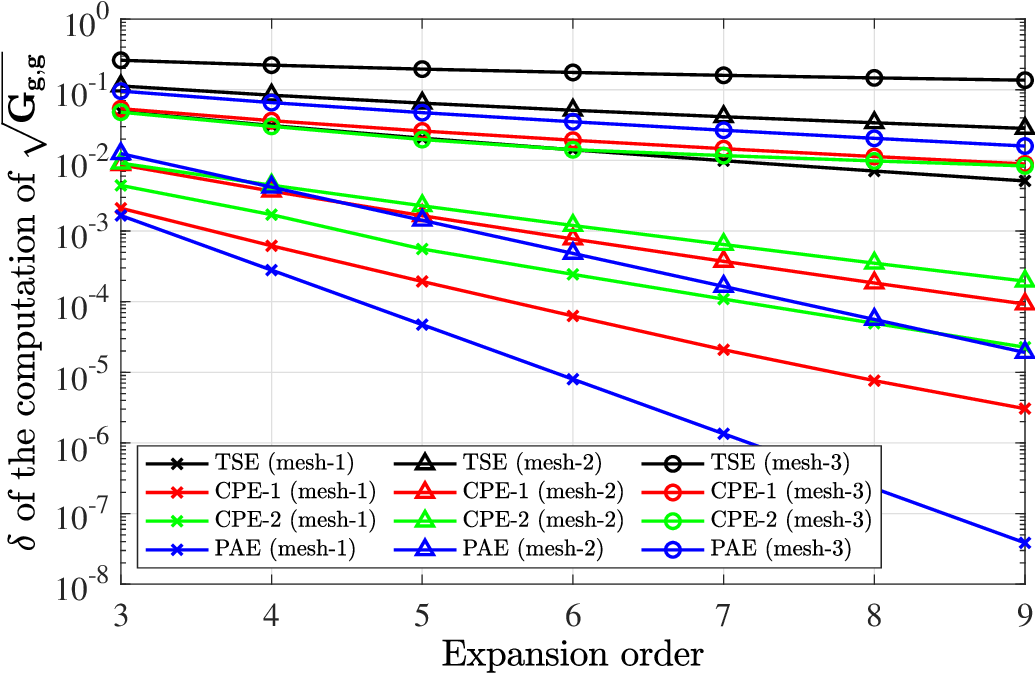}}
\subfigure[]{\includegraphics[width=0.9\columnwidth]{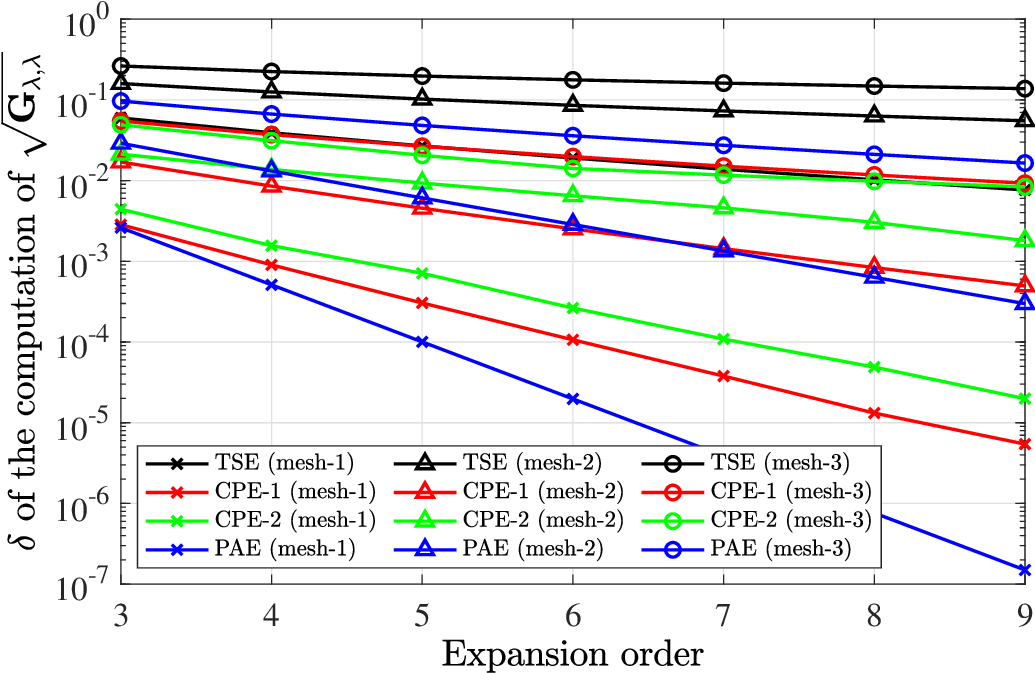}}
\subfigure[]{\includegraphics[width=0.9\columnwidth]{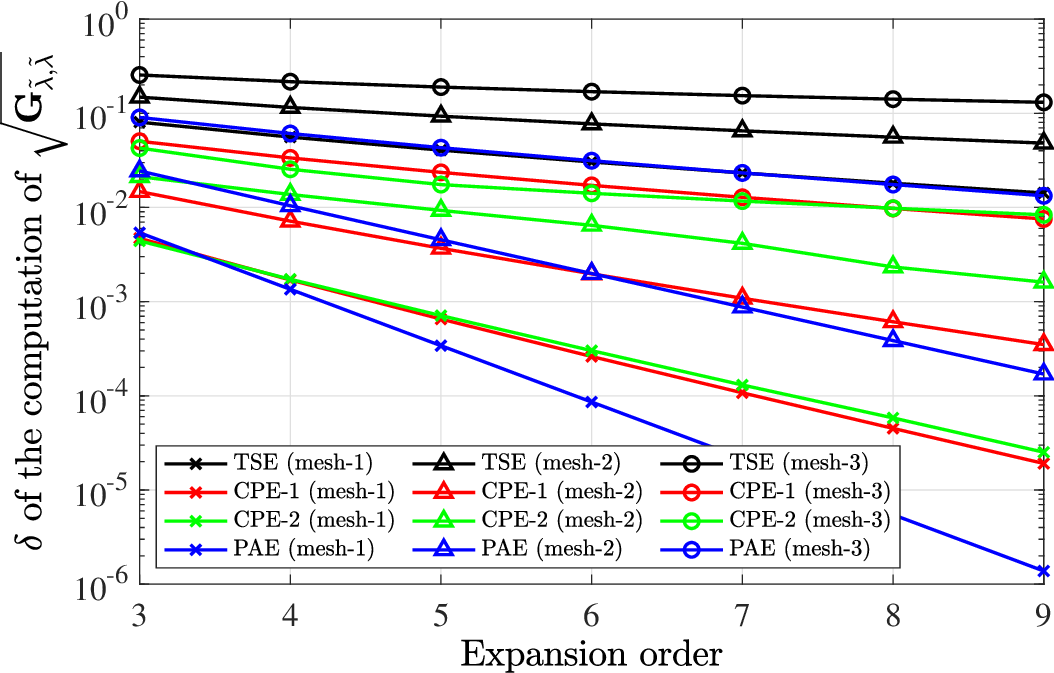}}
\caption{Relative error of the computation of (a) $\sqrt{\mathbf{G}_{\mathbf{f}, \mathbf{f}}}$, (b) $\sqrt{\mathbf{G}_{\mathbf{g}, \mathbf{g}}}$, (c) $\sqrt{\mathbf{G}_{\lambda , \lambda}}$, and (d) $\sqrt{\mathbf{G}_{\tilde{\lambda},\tilde{\lambda}}}$ using TSE, CPE-1, CPE-2, and PAE for mesh-1, mesh-2, and mesh-3 against the expansion order.}
\label{fig1}
\end{figure}

\begin{figure}
\centering
\subfigure[]{\includegraphics[width=0.9\columnwidth]{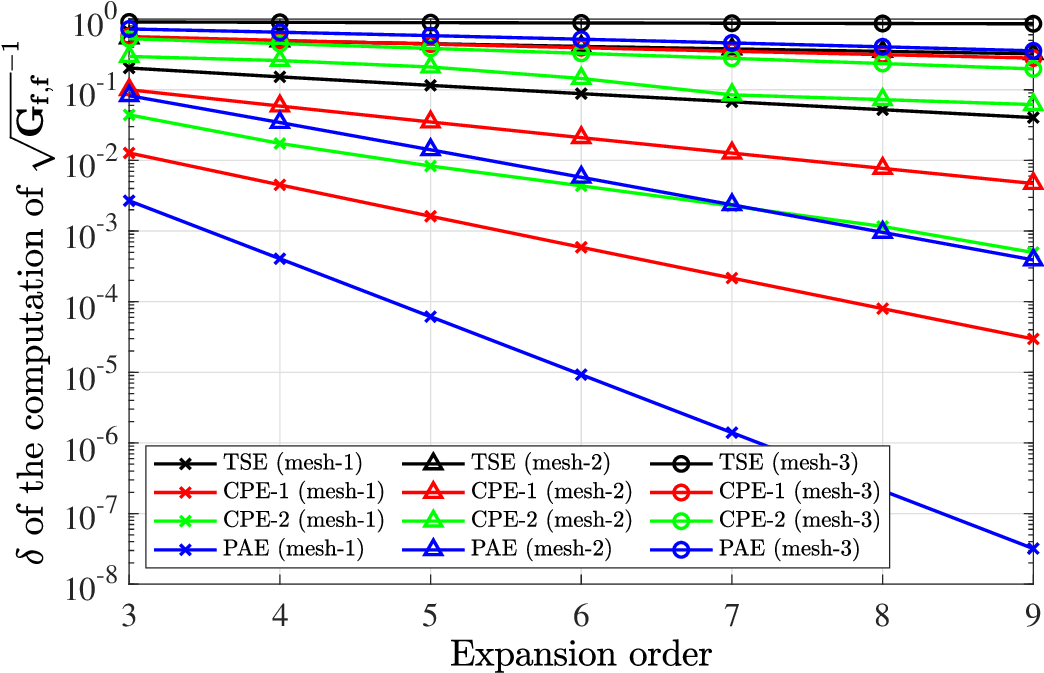}}
\subfigure[]{\includegraphics[width=0.9\columnwidth]{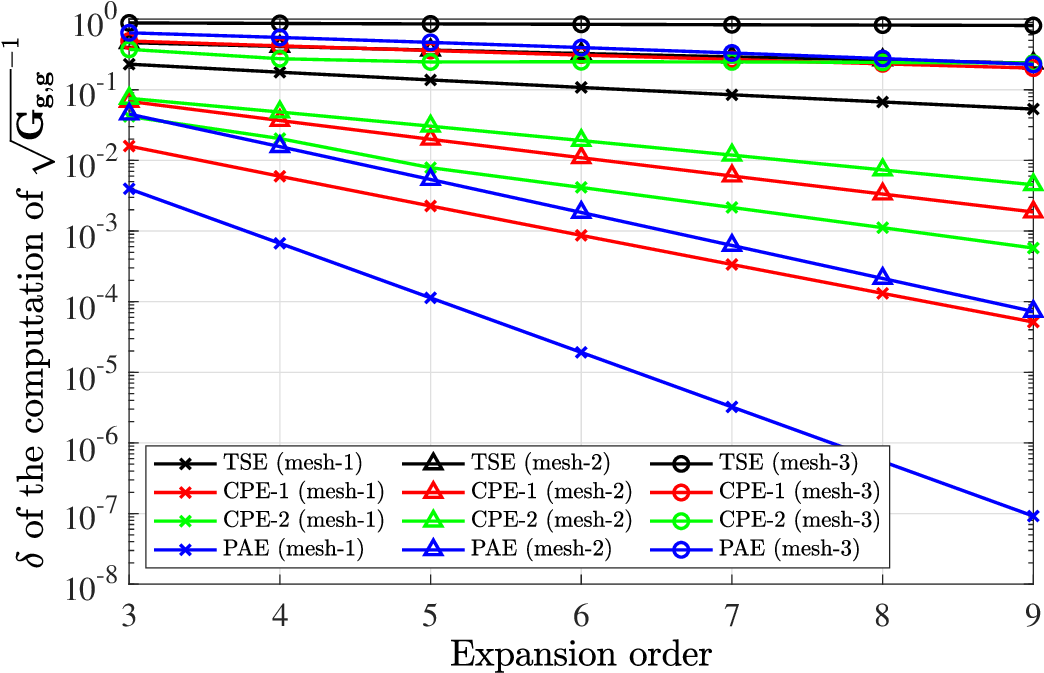}}
\subfigure[]{\includegraphics[width=0.9\columnwidth]{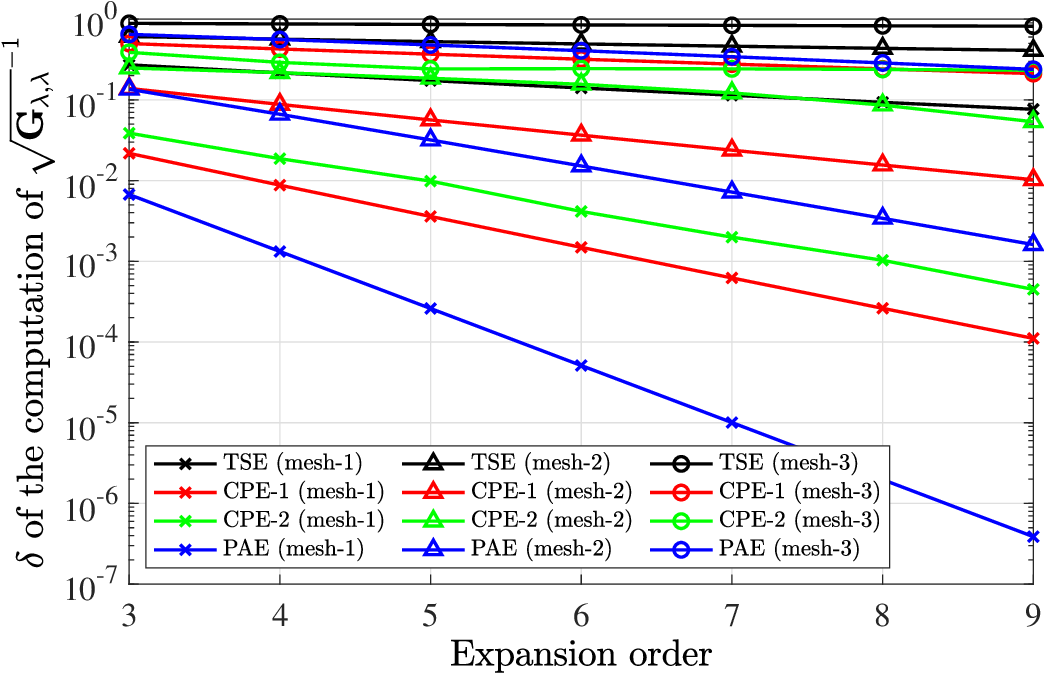}}
\subfigure[]{\includegraphics[width=0.9\columnwidth]{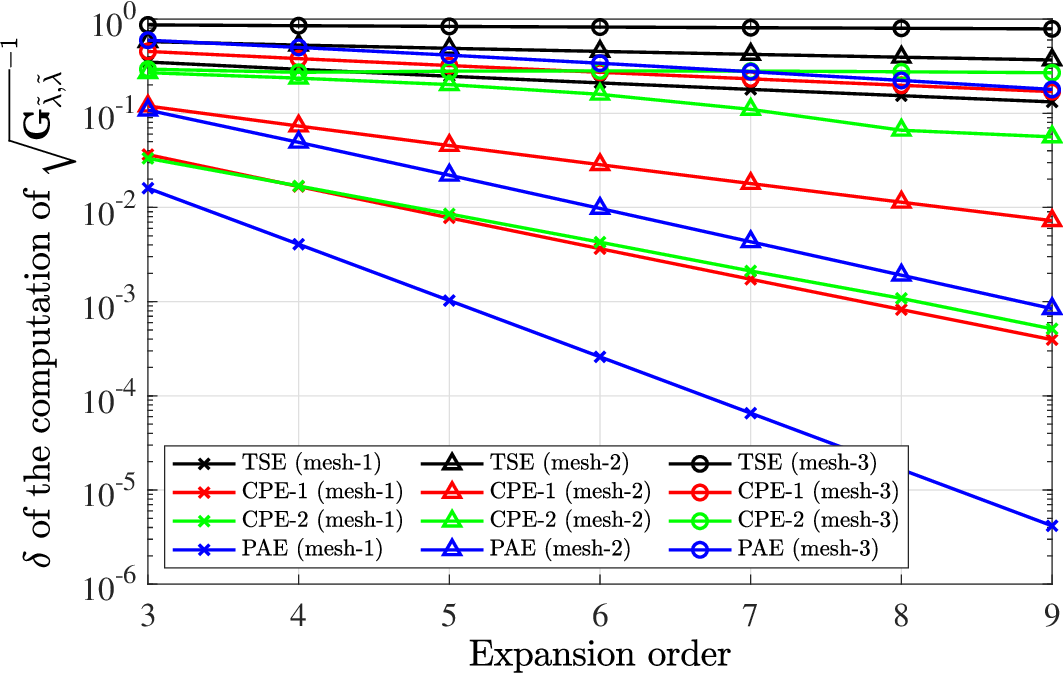}}
\caption{Relative error of the computation of (a) $\sqrt{\mathbf{G}_{\mathbf{f}, \mathbf{f}}}^{-1}$, (b) $\sqrt{\mathbf{G}_{\mathbf{g}, \mathbf{g}}}^{-1}$, (c) $\sqrt{\mathbf{G}_{\lambda , \lambda}}^{-1}$, and (d) $\sqrt{\mathbf{G}_{\tilde{\lambda},\tilde{\lambda}}}^{-1}$ using TSE, CPE-1, CPE-2, and PAE for mesh-1, mesh-2, and mesh-3 against the expansion order.}
\label{fig2}
\end{figure}

As the first numerical experiment, the accuracy of the computation of $\sqrt{\mathbf{G}}$ and $\sqrt{\mathbf{G}}^{-1}$ using the proposed algorithms is investigated. Fig.~\ref{fig1}(a)--(d) and Fig.~\ref{fig2}(a)--(d) show the relative error of the computation of $\sqrt{\mathbf{G}_{\mathbf{f}, \mathbf{f}}}$, $\sqrt{\mathbf{G}_{\mathbf{g}, \mathbf{g}}}$, $\sqrt{\mathbf{G}_{\lambda , \lambda}}$,  $\sqrt{\mathbf{G}_{\tilde{\lambda},\tilde{\lambda}}}$ and $\sqrt{\mathbf{G}_{\mathbf{f}, \mathbf{f}}}^{-1}$, $\sqrt{\mathbf{G}_{\mathbf{g}, \mathbf{g}}}^{-1}$, $\sqrt{\mathbf{G}_{\lambda , \lambda}}^{-1}$,  $\sqrt{\mathbf{G}_{\tilde{\lambda},\tilde{\lambda}}}^{-1}$ for the three sets of meshes against the expansion order, respectively. As expected, the accuracy of the computation of $\sqrt{\mathbf{G}}$ and $\sqrt{\mathbf{G}}^{-1}$ using the proposed algorithms for all the sets of meshes increases with the increase of the expansion order. In addition, for a fixed expansion order and a given algorithm, the accuracy of the computation of $\sqrt{\mathbf{G}}$ and $\sqrt{\mathbf{G}}^{-1}$ for mesh-1 is higher than that for mesh-2 and mesh-3. This effect is a consequence of the lower condition numbers of the Gram matrices for mesh-1  than that for the other two sets of meshes.
Moreover, for the low condition numbers (e.g., mesh-1), PAE converges faster than TSE and CPE. 
Fig.~\ref{fig1}(a)--(d) and Fig.~\ref{fig2}(a)--(d) also show that, for the same Gram matrix, the accuracy of the computation of $\sqrt{\mathbf{G}}^{-1}$ is generally lower than that of  $\sqrt{\mathbf{G}}$ using the same expansion order for a given algorithm and a given set of meshes---in accordance with the report in Table \ref{tab:order_chebyshev_sqrt}. 
In addition, the figures show that we could use CPE-2 as a good approximation of CPE-1 by directly using the bounds of the living known ranges of $n_0$ and their corresponding coefficients in Table~\ref{tab:coefs_chebyshev_sqrt} and Table~\ref{tab:coefs_chebyshev_isqrt} for CPE for convenience. The figures also show that the error behaviour of CPE-1 and CPE-2 converges towards one another when the actual value of $n_0$ for a given set of meshes and a given Gram matrix corresponds to the value on $n_0$ used to select the CPE-2 tabulated coefficients.
\color{black}

\subsection{Application to Spectral Analysis}

As the second numerical experiment, the use of the proposed technique as a tool for the study and manipulation of the spectrum of the relevant integral operators 
is showcased for the problem of scattering by a PEC sphere of radius \SI{1}{\meter}. The surface of the PEC sphere is uniformly discretized with $\{ N_{\mathrm{E}}, N_{\mathrm{P}}, N_{\mathrm{V}}\}=\{ \num{1080}, \num{719}, \num{361}  \}$ and the wavenumber is set to $k = \SI{0.1}{\radian\per\meter}$. 

\begin{figure}
\centering
\includegraphics[width=0.853\columnwidth]{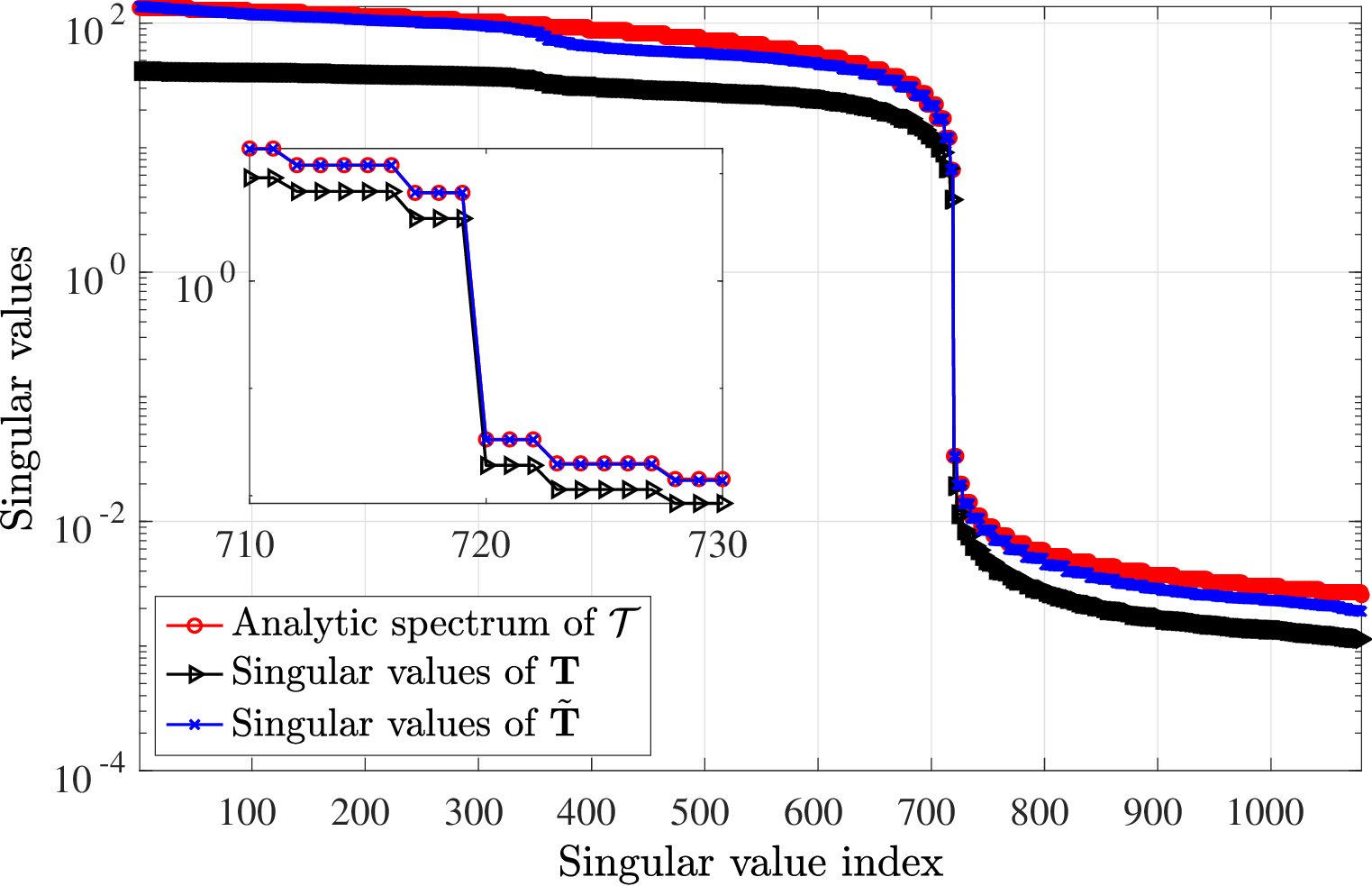}
\caption{Comparison of the singular values of $\mathbf{T}$ and $\tilde{\mathbf{T}}$ with the associated analytic spectrum of $\mathcal{T}$ including a zoom around singular value index \num{720}.}
\label{fig3}
\end{figure}

Fig.~\ref{fig3} compares both the singular values of the standard EFIE boundary element matrix $\mathbf{T}$ and its normalized counterpart $\tilde{\mathbf{T}}=\sqrt{\mathbf{G}_{\mathbf{f}, \mathbf{f}}}^{-1} \mathbf{T} \sqrt{\mathbf{G}_{\mathbf{f}, \mathbf{f}}}^{-1}$ with the analytic spectrum of $\mathcal{T}$ obtained from a spherical harmonics analysis of the problem.
The results are obtained with PAE for the numerical approximation of $\sqrt{\mathbf{G}_{\mathbf{f}, \mathbf{f}}}^{-1}$ with $N_\mathrm{A}=9$.
As seen in Fig.~\ref{fig3}, the  singular values of the normalized  matrix $\tilde{\mathbf{T}}$ show a good agreement with the analytic spectrum of $\mathcal{T}$ while that of the standard matrix $\mathbf{T}$ do not.

\section{Conclusion}
\label{sec:conclusion}
In this work, three algorithmic approaches are delineated for the computation of the square roots and reciprocal square roots of four relevant Gram matrices used for the zeroth-order discretization of the boundary element matrices in electromagnetics. 
General algorithms have been tailored to the specific integral equation scenario and tables containing the different sets of the expansion coefficients have been 
provided along with the comparative numerical experiments that evidence
the advantages and disadvantages of the different strategies.

\end{document}